
\input amstex
\documentstyle{amsppt}
\magnification=1200
 \vsize19.5cm
  \hsize13.5cm
   \TagsOnLeft
\pageno=1
\baselineskip=15.0pt
\parskip=3pt

\def\p{\partial}
\def\noo{\noindent}
\def\eps{\varepsilon}
\def\lam{\lambda}
\def\Om{\Omega}

\def\pom{{\p \Om}}
\def\bom{{\overline\Om}}
\def\R{\Bbb R}

\def\A{ Y}
\def\B{\psi}

\def\wtt{\widetilde}

\def\dist{\text{dist}}
\def\det{\text{det}}
\def\log{\text{log}}

\def\ol{\overline}

\def\D{\nabla}
\def\phi{\varphi}

\nologo \NoRunningHeads

\topmatter

\title{On the second boundary value problem for Monge-Amp\`ere
type equations and optimal transportation}\endtitle
\author{Neil S. Trudinger\ \ \ \
       Xu-Jia Wang}\endauthor

\affil{Centre for Mathematics and Its Applications\\
       The Australian National University\\
       Canberra, ACT 0200\\
       Australia\\ }\endaffil

\address{Centre for Mathematics and Its Applications,
       Australian National University,
       Canberra, ACT 0200,
       Australia}\endaddress

\noindent \email{\newline Neil.Trudinger\@anu.edu.au;\ \ \
wang\@maths.anu.edu.au}\endemail

\thanks{This research  was supported
        by the Australian Research Council.}\endthanks

\abstract { This paper is concerned with the existence of globally
smooth solutions for the second boundary value problem for
Monge-Amp\`{e}re type equations and the application to regularity
of potentials in optimal transportation. The cost functions
satisfy a weak form of the condition A3, which was introduced
 in a
recent paper with Xi-nan Ma  in conjunction with interior
regularity. Consequently they include the quadratic cost function
case of Caffarelli and Urbas as well
 as the various examples in the earlier
work. The approach is through the derivation of global estimates for
second derivatives of solutions. }\endabstract

\endtopmatter

\vskip-20pt

\document

\baselineskip=14.2pt
\parskip=3pt

\centerline {\bf  1. Introduction}

\vskip10pt

This paper is concerned with the global regularity of solutions of
the second boundary value problem for equations of Monge-Amp\`ere
type and its application to the regularity of potentials in
optimal transportation problems with non-quadratic cost functions.

The Monge-Amp\`ere equations under consideration have the general
form
$$\det \{ D^2 u - A(\cdot, u, Du)\}=B(\cdot, u, Du), \tag 1.1$$
where $A$ and $B$ are given $n\times n$ matrix and scalar valued
function defined on $\Om\times \R\times \R^n$, where $\Om$ is a
domain in Euclidean $n$-space, $\R^n$. We use $(x, z, p)$ to
denote points in $\Om\times \R\times \R^n$ so that $A(x, z, p)\in
\R^n\times \R^n$, $B(x, z, p)\in \R$ and $(x, z, p)\in \Om\times
\R\times \R^n$. The equation (1.1) will be elliptic, (degenerate
elliptic), with respect to a solution $u\in C^2(\Om)$ whenever
$$D^2 u-A(\cdot, u, Du)>0\ \ \ (\ge 0), \tag 1.2$$
whence also $B>0$\ $(\ge 0)$.

A particular form  of (1.1) arises from the prescription of the
Jacobian determinant of a mapping $Tu$ defined by
$$T u=\A(\cdot, u, Du), \tag 1.3$$
where $\A$ is a given vector valued function on $\Om\times
\R\times \R^n$, namely
$$\det\, D\A(\cdot, u, Du)=\B(\cdot, u, Du). \tag 1.4$$
Assuming that the matrix
$$\A_p= [D_{p_j}\A^i]\tag 1.5$$
is  non-singular, we may write (1.4) in the form (1.1), that is,
$$\det \{D^2 u+\A^{-1}_p(\A_x+\A_z\otimes Du)\}
             =\frac {\B}{|\det\, \A_p|}, \tag 1.6$$
for degenerate elliptic solutions $u$.

The {\it second boundary value problem} for equation (1.4) is to
prescribe the image
$$T u(\Om)=\Om^*,\tag 1.7$$
where $\Om^*$ is a given domain in $\R^n$. When $\A$ and $\B$ are
independent of $z$ and $\B$ is separable in the sense that
$$\B(x, p)=f(x)/g \circ \A(x, p) \tag 1.8$$
for positive $f, g\in L^1(\Om)$, $L^1(\Om^*)$ respectively, then a
necessary condition for the existence of an elliptic solution,
for which the mapping $T$ is a diffeomorphism, to
the second boundary value problem (1.4) (1.7) is the {\it mass
balance} condition
$$\int_\Om f=\int_{\Om^*} g. \tag 1.9$$
The second boundary value problem (1.4) (1.7) arises naturally in
optimal transportation. Here we are given a cost function $c:\
\R^n\times \R^n\to \R$ and the vector field $\A$ is generated by
the equation
$$c_x(x, \A(x, p))=  p, \tag 1.10$$
which we assume to be uniquely solvable for $p\in\R^n$, with
non-vanishing determinant, that is
$$\det \, c_{x, y}(x, y)\ne 0\tag 1.11$$
for all $x, y\in \Om\times \Om^*$. Using the notation
$$c_{i j\cdots, kl\cdots}=\frac{\p}{\p x_i}\frac{\p}{\p x_j}\cdots
 \frac{\p}{\p y_k}\frac{\p}{\p y_l}\cdots c\tag 1.12$$
we have
$$\A_p(x,p)= [c^{i, j}(x, \A(x, p))],\tag 1.13$$
where $[c^{i, j}]$ is the inverse of $[c_{i, j}]$. The
corresponding Monge-Amp\`ere equation can now be written as
$$\det \{D^2 u-c_{xx}(\cdot, \A(\cdot, Du))\}
                       =|\det \, c_{x, y}| \B,     \tag 1.14$$
that is in the form (1.1) with
$$\align
A(x, z, p) & =c_{xx}(x, \A(x, p)), \tag 1.15\\
B(x, z, p) & =|\det\,  c_{x, y}(x, \A(x, p))|\B(x, z, p).\\
\endalign $$
In the case of the (quadratic) cost function
$$c(x, y)=   x\cdot y, \tag 1.16$$
we have
$$\A(x, p)=p,\ \ \ \ T u=Du, \tag 1.17$$
and equation (1.14) reduces to the standard Monge-Amp\`ere
equation
$$\det\, D^2 u=\B. \tag 1.18$$
For this case global regularity of solutions was proved by
Delano\"e [D], Caffarelli [C2] and Urbas [U1], with (conditional)
interior regularity shown earlier by Caffarelli [C1]. In this
paper we will prove global estimates and regularity under
corresponding conditions. In particular, we will assume that the
cost function $c\in C^4(\R^n\times \R^n)$ satisfies the following
conditions:

\noo (A1)\ \  For each $p, q\in \R^n$, there exists unique
$y=\A(x, p)$, $x=X(q, y)$ such that
$$\align
c_x(x, y) & =  p\ \ \ \forall \ x\in\Om,\\
c_y(x, y) & =  q\ \ \ \forall\ y\in \Om^*.\\
\endalign $$
$$\det\, c_{x, y}(x, y)\ne 0,\ \ \
        \forall\ x\in\bom, y\in \bom^*. \tag A2$$
$$\align
\Cal F(x,p;\xi,\eta): & = D_{p_ip_j}A_{kl}(x, p)\xi_i\xi_j\eta_k\eta_l \tag A3w\\
& \ge 0\ \ \
   \forall\ \ x\in \Om, p\in \R^n, \ \xi\perp \eta\in\R^n.\\
   \endalign $$

Conditions A1 and A2 are the same conditions as in [MTW]
but condition A3w is the degenerate form of condition A3 in [MTW].
As will be seen in our examples, we do not necessarily require $c$
to be defined on all of $\R^n\times \R^n$ and the vectors $p$ and $q$
in conditions A1 and A3w need only lie in the ranges of $c_x(x, y)$ and
$c_y(x, y)$ on $\Om\times\Om^*$. Moreover, as done at the 
outset in [MTW], we may also write
$$\Cal F(x,p;\xi,\eta) = (c_{ij,rs} - c^{k^\prime,l^\prime}
c_{ij,k^\prime}c_{l^\prime,rs})c^{r,k}c^{s,l}(x,y)\xi_i\xi_j\eta_k\eta_l\tag 1.19$$
where $y$ and $p$ are related through A1. This shows
that condition A3w is also symmetric in x and y.

In our paper [MTW], we also introduced a notion of convexity of
domains with respect to cost functions, namely $\Om$ is
$c$-convex, with respect to $\Om^*$, if the image $c_y(\cdot,
y)(\Om)$ is convex in $\R^n$ for each $y\in \Om^*$, while
analogously $\Om^*$ is $c^*$-convex, with respect to $\Om$, if the
image $c_x(x,\cdot)(\Om^*)$ is convex for each $x$ in $\Om$. For
global regularity we need to strengthen these conditions in the
same way that convexity is strengthened to uniform convexity.
Namely we define $\Om$ to be uniformly $c$-convex, with respect to
$\Om^*$, if $\Om$ is $c$-convex, with respect to $\Om^*$, $\pom\in
C^2$ and there exists a positive constant $\delta_0$ such that
$$[D_i\gamma_j(x)-c^{l, k} c_{ij, l}(x, y) \gamma_k(x)]
                           \tau_i\tau_j(x) \ge \delta_0\tag 1.20$$
for all $x\in\pom$, $y\in\Om^*$, unit tangent vector $\tau$ and
outer unit normal $\gamma$. By pulling back with the mappings
$c_y(\cdot , y)$, we see that this is equivalent to the condition
that the image domains $c_y(\cdot, y)(\Omega)$ be uniformly convex
with respect to $y \in \Omega^*$. Similarly we call $\Om^*$
uniformly $c^*$-convex, with respect to $\Om$, when $c^*(x,
y)=c(y, x)$. Note that if $\Omega$ is simply connected with
boundary $\partial \Omega \in C^2$, then $\Omega$ is $c-$convex if
and only if (1.20) holds for $\delta_0=0$.

We can now formulate our main estimate.

\proclaim{Theorem 1.1} Let $c$ be a cost function satisfying
hypotheses A1, A2, A3w, with respect to bounded $C^4$ domains
$\Om, \Om^*\in \R^n$ which are respectively uniformly $c$-convex,
$c^*$-convex with respect to each other. Let $\B$ be a positive
function in $C^2(\bom\times \R\times \R^n)$. Then any elliptic
solution $u\in C^3(\bom)$ of the second boundary value problem
(1.14), (1.7) satisfies the a priori estimate
$$ |D^2 u|\le C, \tag 1.21$$
where $C$ depends on $c,\B,\Om$ , $\Om^*$ and $\sup_\Om|u|$.
\endproclaim

As we will indicate later, the smoothness assumption on the
solution and the data may be reduced. Further regularity also
follows from the theory of linear elliptic equations for example
if $c, \Om, \Om^*, \B$ are $C^\infty$ then the solution $u\in
C^\infty(\bom)$.The dependence of the estimate (1.21) on
$\sup_\Om|u|$ may be removed if $\B$ is independent of $u$
as in (1.8).

As a consequence of Theorem 1.1, we may conclude existence
 theorems for classical solutions by the method of continuity. 

\proclaim {Theorem 1.2} Suppose in addition to the above
hypotheses that the function $\B$ satisfies (1.8) (1.9). Then
there exists a unique (up to additive constants) elliptic solution
$u\in C^3(\bom)$ of the second boundary value problem (1.14),
(1.7).
\endproclaim

From Theorem 1.2, we also obtain an existence result for classical
solutions of the Monge-Kantorovich problem in optimal
transportation. As above we let $c\in C^4(\R^n\times \R^n)$ be a
cost function and $\Om, \Om^*$ be two bounded domains in $\R^n$
satisfying the hypotheses of Theorem 1.1. Let $f>0, \in
C^2(\bom)$, $g>0, \in C^2(\bom^*)$ be positive densities
satisfying the mass balance condition (1.9). Then the
corresponding optimal transportation problem is to find a measure
preserving mapping $T_0:\ \Om\to\Om^*$ which maximizes the cost
functional
$$C(T)= \int_\Om f(x) c(x, T(x)) dx\tag 1.22$$
among all measure preserving mappings $T$ from $\Om$ to $\Om^*$. A
mapping $T:\ \Om\to\Om^*$ is called measure preserving if it is
Borel measurable and for any Borel set $E\subset \Om^*$,
$$\int_{T^{-1}(E)} f=\int_E g.\tag 1.23$$

\proclaim {Theorem 1.3} Under the above hypotheses, there exists a
unique diffeomorphism $T\in [C^2(\bom)]^n$ maximizing the
functional (1.22), given by
$$T(x) = \A(x, Du(x)), \tag 1.24$$
where $u$ is an elliptic solution of the boundary problem (1.7),
(1.14).
\endproclaim

The solution $u$ of (1.7), (1.14) is called a potential. Note that in [MTW]
and elsewhere the cost functions and potentials are the negatives of those
here and the optimal transportation problem is written, (in its usual form),
as a minimization problem.

The plan of this paper is as follows. In Section 2, we prove that
boundary conditions of the form (1.7) are oblique with respect to
functions for which the Jacobian $DT$ is non-singular and we
estimate the obliqueness for solutions of the boundary value
problem (1.14), (1.7) under hypotheses A1 and A2, (Theorem 2.1).
Here the twin assumptions of $\Om$ and $\Om^*$ being uniformly
 $c$ and $c^*$-convex with respect to each other are critical. In
Section 3, we prove that second derivatives of solutions of
equation (1.14) can be estimated in terms of their boundary values
under hypothesis A3w, (Theorem 3.1). This estimation is already
immediate from [MTW] when the non-degenerate condition A3 is
satisfied. The argument is carried out for equations of the
general form (1.1) (with symmetric $A$), in the presence of a
global barrier, which is not necessary
 in the optimal
transportation case (Theorem 3.2). This estimation also arises in
the treatment of the classical Dirichlet problem [T3]. The
proof of the global second derivative estimates in Theorem 1.1
for solutions of the boundary value problem (1.14), (1.7)
 is completed in Section 4. Here the procedure is similar to that 
 in [LTU] and [U1].We remark here that this global estimate also extends to
 the more  general prescribed Jacobian equation (1.6) [T3]. In
Section 5, we commence the proof of the existence result, Theorem
1.2, by adapting the method of continuity [GT] and establish the result
under a stronger uniform c-convexity hypothesis, 
(employed in earlier versions of this paper).
 Section 6 is devoted to the applications to optimal
transportation and the proof of Theorem 1.3 from Theorem 1.2, 
which implies the
global regularity of the potential functions in [MTW], under
condition A3w.  In Section 7, we  finally complete the proof of
Theorem 1.2  in its full generality, by showing that there exists
a smooth function satisfying the ellipticity condition (1.2), together
with the boundary condition (1.7), (at least for approximating domains).
In the last section, we discuss our results in the light of examples, 
most of which are already given in [MTW].

 \vskip20pt

\centerline{\bf 2. Obliqueness}

\vskip10pt

In this section, we prove that the boundary condition (1.7)
implies an oblique boundary condition and estimate the
obliqueness. First we recall that a boundary condition of the form
$$G(x, u, Du)=0\ \ \ \text{on}\ \ \pom\tag 2.1$$
for a second order partial differential equation in a domain $\Om$
is called oblique if
$$G_p\cdot\gamma >0\tag 2.2$$
for all $(x, z, p)\in\pom\times\R\times \R^n$, where $\gamma$
denotes the unit outer normal to $\pom$. Let us now assume that
$\phi$ and $\phi^*$ are $C^2$ defining functions for $\Om$ and
$\Om^*$ respectively, with $\phi, \phi^*<0$ near $\pom, \pom^*$,
$\phi=0$ on $\pom$, $\phi^*=0$ on $\pom^*$, $\D \phi, \D \phi^*\ne
0$ near $\pom, \pom^*$. Then if $u\in C^2(\bom)$ is an elliptic
solution of the second boundary value problem (1.4), (1.7), we
must have
$$\phi^*(T u)=0\ \ \ \text{on}\ \ \pom,\ \ \
   \phi^*(T u)<0\ \ \ \text{near}\ \ \pom. \tag 2.2$$
By tangential differentiation, we obtain
$$\phi_i^*(D_jT^iu)\tau_j=0\tag 2.3$$
for all unit tangent vectors $\tau$, whence
$$\phi^*_i(D_jT^i)=\chi\gamma_j\tag 2.4$$
for some $\chi\ge 0$. Consequently
$$\phi_i^*c^{i, k}(u_{jk}-c_{jk})=\chi\gamma_j, \tag 2.5$$
that is
$$\phi^*_ic^{i, k} w_{jk}=\chi\gamma_j, \tag 2.6$$
where
$$w_{ij}=u_{ij}-c_{ij}. \tag 2.7$$
At this point we observe that $\chi>0$ on $\pom$ since $|\D
\phi^*|\ne 0$ on $\pom$ and $\det\,DT \ne 0$. Using the
ellipticity of (1.4) and letting $[w^{ij}]$ denote the inverse
matrix of $[w_{ij}]$, we then have
$$\phi_i^*c^{i, k}=\chi w^{jk}\gamma_j. \tag 2.8$$
Now writing
$$G(x, p)=\phi^*(\A(x, p)),\tag 2.9$$
we have
$$\beta_k:=G_{p_k}(\cdot, Du)=\chi w^{jk}\gamma_j, \tag 2.10$$
whence
$$\align
\beta\cdot \gamma & =\chi w^{ij}\gamma_i\gamma_j\tag 2.11\\
                  & > 0\\
                  \endalign $$
on $\pom$. We obtain a further formula for $\beta\cdot\gamma$,
from (2.6), namely
$$\align
\phi^*_ic^{i, k}w_{jk}\phi^*_lc^{l, j}
                  &=\chi\phi_l^*c^{l,j}\gamma_j\tag 2.12\\
                  &=\chi(\beta\cdot \gamma).\\
                  \endalign $$
Eliminating $\chi$ from (2.11) and (2.12), we have
$$(\beta\cdot\gamma)^2=(w^{ij}\gamma_i\gamma_j)
             (w_{kl}c^{i,k}c^{j, l} \phi_i^*\phi_j^*) .\tag 2.13$$
We call (2.13) a formula of Urbas type, as it was proved by Urbas
[U1] for the special case, $c(x, y)=x\cdot y$, $\A(\xi, p)=p$, of
the Monge-Amp\`ere equation. Note that to prove (2.13), we only
used conditions A1 and A2 and moreover (2.13) continues to hold in
the generality of (1.4).

Our main task now is to estimate $\beta \cdot \gamma$ from below
for solutions of (1.14), (1.7). For this in addition to conditions
A1, A2, we also need the uniform $c$ and $c^*$ convexity of $\Om$
and $\Om^*$ respectively. Our approach is similar to [U1] for the
special case of the Monge-Amp\`ere equation and begins by invoking
the key idea from [T] for estimating double normal derivatives of
solution of the Dirichlet problem. Namely we fix a point $x_0$ on
$\pom$ where $\beta\cdot\gamma$ is minimized, for an elliptic
solution $u\in C^3(\bom)$, and use a comparison argument to
estimate $\gamma\cdot D(\beta\cdot\gamma)$ from above. Without
some concavity condition in $p$ the quantity $\beta\cdot \gamma$
does not satisfy a nice differential inequality so we will get
around this by considering instead the function
$$v=\beta\cdot \gamma-\kappa\phi^*(T u)\tag 2.14$$
for sufficiently large $\kappa$, where now the defining function
$\phi^*$ is chosen so that
$$(D_{ij}\phi^*(Tu)-c^{k,l}c_{l, ij}(\cdot, Tu)D_k\phi^*(Tu)\xi_i\xi_j
  \ge \delta_0^*|\xi|^2\tag 2.15$$
near $\pom$, for all $\xi\in \R^n$ and some positive constant
$\delta_0^*$. Inequality (2.15) is possible by virtue of the
uniform $c^*$-convexity of $\Om^*$, with the function $\phi^*$
given, for example by
$$\phi^*=-ad^*+b(d^*)^2,\tag 2.16$$
where a and b are positive constants and
 $d^*$ denotes the distance function for $\Om^*$, [GT].

By differentiation of equation
(1.14), in the form (1.1), we obtain, for $r=1,\cdots, n$,
$$w^{ij}\{D_{ij}u_r-D_{p_k}A_{ij}(x, Du) D_ku_r
             -D_{x_r}A_{ij}(x, Du)\}=D_r\log B. \tag 2.17$$
Introducing the linearized operator $L$,
$$Lv=w^{ij} (D_{ij}v-D_{p_k}A_{ij}D_kv), \tag 2.18$$
we need to compute $Lv$ for $v$ given by (2.14). Setting
$$F(x, p)=G_p(x, p)\cdot \gamma(x)-\kappa G(x, p), \tag 2.19$$
where $G$ is defined by (2.9), we see that
$$v(x)=F(x, Du(x)).\tag 2.20$$
Writing $b^{ij}_k= -D_{p_k}A_{ij}$, we then have
$$\align
Lv= & w^{ij}\{F_{p_r}D_{ij}u_r+F_{p_rp_s}D_{ir}u D_{js}u\\
    & +F_{x_ix_j}+2F_{x_ip_r}D_{jr}u
    +b^{ij}_k(F_{x_k}+F_{p_r}D_ku_r)\}\tag 2.21\\
    \endalign$$
In the ensuing calculations, we will often employ the following
formulae,
$$\align
c^{i,j}_k(x, y) &=D_{x_k}c^{i, j}(x, y)\\
          &= -c^{i, l}c^{r, j}c_{kl, r}(x, y),\tag 2.22\\
c^{i,j}_{,k}(x, y) &=D_{y_k}c^{i, j}(x, y)\\
          &= -c^{i, l}c^{r, j}c_{l, kr}(x, y),\\
          \endalign $$
as well as (1.13). Indeed, using (1.13) and (2.22), we have
$$\align
G_{p_ip_j}=&D_{p_j}(\phi^*_kc^{k, i})\tag 2.23\\
  =& \phi^*_{kl}c^{k, i}c^{l, j}-\phi^*_kc^{s,j}c^{k, i}_{,s}\\
  =& c^{k, i}c^{l, j}\{\phi^*_{kl}-\phi^*_rc^{r, s}c_{s, kl}\}\\
  \endalign$$
so that
$$\align
G_{p_ip_j}(x, Du) \xi_i\xi_j
 \ge & \delta_0^* \sum |c^{i,j} \xi_j|^2\tag 2.24\\
  \ge & \kappa_0^*\, |\xi|^2\\
  \endalign $$
for a further positive constant $\kappa_0^*$. By choosing $\kappa$
sufficiently large, we can then ensure that
$$F_{p_ip_j}(x, Du)\xi_i\xi_j\le -\frac 12 \kappa |\xi|^2\tag 2.25$$
near $\pom$. Substituting into (2.20) and using (2.16), it follows
that
$$Lv\le -\frac 14 \kappa w_{ii}+C(w^{ii}+1)+D_{p_k} log B D_kv, \tag 2.26$$
where $C$ is a constant depending on $c, \B, \Om$ and $\Om^*$, as
well as $\kappa$.

Next we observe that unless the defining function $\phi^*$ extends to
all of $\Om^*$ so that  (2.15) is satisfied for all $Tu\in\Om^*$, we have no 
control on the neighbourhood of $\pom$, where (2.26) holds. This is
remedied by replacing $G$ in (2.19) by a function satisfying (2.24) in
all of $\Om$ , agreeing with (2.9) near $\pom$, for example by taking
 $$G(x,p) = m_h\{(\phi^*(\A(x, p)),a_1(|p|^2- K^2)\}, \tag 2.27$$
where $a_1$ and $K$ are positive constants, with $a_1$ sufficiently
small and $K > $max$|Du|$, and for $h$ sufficiently small,  $m_h$ is the 
mollification of the max-function of two variables.

A suitable barrier is now provided by the uniform $c$-convexity of
$\Om$ which implies, analogously to the case of $\Om^*$ above,
that there exists a defining function $\phi$ for $\Om$ satisfying
$$[D_{ij}\phi-c^{l, k}c_{ij, l}(x, Tu)D_k\phi]\xi_i\xi_j
       \ge \delta_0|\xi|^2, \tag 2.28$$
in a fixed neighbourhood of  $\pom$, (for a constant $\delta_0>0$).
 By appropriate choice, of say the constants $a$ and $b$ in (2.16),
 without the $*$, (or following the uniformly convex case 
 in [GT, Chapter 14]),
 we may obtain, by virtue of (2.21),
$$L\phi\ge \delta_0w^{ii}+Kw^{ij}D_i\phi D_j\phi,\tag 2.29$$
for a given constant K.
Combining (2.26) and (2.29), and using the positivity of $B$, we
then infer by the usual barrier argument,( which entails fixing a 
small enough neighbourhood of $\pom$, [GT]),
$$\gamma\cdot Dv(x_0)\le C, \tag 2.30$$
where again $C$ is a constant depending on $c, \Om, \Om^*$ and
$\B$. From (2.28) and since $x_0$ is a minimum point of $v$ on
$\pom$, we can write
$$Dv(x_0)=\tau \gamma(x_0)\tag 2.31$$
where $\tau\le C$. To use the information embodied in (2.31), we
need to calculate
$$\align
D_i(\beta\cdot\gamma)
   =& D_i\{\phi_k^* c^{k,j}\gamma_j\}\tag 2.32\\
   =& \phi^*_{kl}D_i(T^lu)c^{k, j}\gamma_j
       +\phi^*_kc^{k, j} D_i\gamma_j
          +\phi_k^*\gamma_j(c^{k, j}_i+c^{k, j}_{, l}D_iT^lu)\\
   =& \phi^*_kc^{k, j}(D_i\gamma_j-c^{s,r}c_{ij, s}\gamma_r)
       +(\phi^*_{kl}-\phi^*_rc^{r, s}c_{s,kl})
                     c^{k,j}\gamma_jD_iT^lu\\
       \endalign $$
Multiplying by $\phi^*_t c^{t,i}$ and summing over $i$, we obtain
$$\align
\phi^*_tc^{t, i}D_i(\beta\cdot\gamma)
 =& \phi^*_k\phi^*_t c^{k,j}c^{t,i}
                 (D_i\gamma_j-c^{s, r}c_{ij,s}\gamma_r)\tag 2.33\\
 &+\phi^*_tc^{t,i}c^{k,j}\gamma_jc^{l, m}w_{im}
                 (\phi^*_{kl}-\phi^*_rc^{r,s}c_{s, kl})\\
 \ge & \delta_0\sum |\phi^*_ic^{i, j}|^2\\
 \endalign $$
by virtue of the uniform $c$-convexity of $\Om$, the
$c^*$-convexity of $\Om^*$ and (2.6). Consequently, from (2.19)
and (2.31), we obtain at $x_0$,
$$-\kappa w_{kl}c^{i,k}c^{j,l}\phi^*_i\phi^*_j
 \le C(\beta\cdot\gamma)-\tau_0\tag 2.34$$
for positive constants, $C$ and $\tau_0$. Hence if
$\beta\cdot\gamma\le \tau_0/2C$, we have the lower bound
$$w_{kl}c^{i, k}c^{j, l}\phi^*_i\phi^*_j
                       \ge \frac{\tau_0}{2\kappa}. \tag 2.35$$

\vskip10pt

To complete the estimation of $\beta\cdot\gamma$ we invoke the
dual problem to estimate $w^{ij}\gamma_i\gamma_j$ at $x_0$.
Assuming for the moment that $Tu$ is one to one, we let $u^*$
denote the $c$-transform of $u$, defined for $y = Tu(x)\in\Om^*$
by
$$u^*(y)=c(x, y)-u(x).\tag 2.36$$
It follows that
$$\align
Du^*(y) & =c_y(x, y)\tag 2.37\\
       & =c_y(T^*u^*(y),y), \\
       \endalign $$
where
$$\align
T^*u^*(y) & =X (Du^*, y)\tag 2.38\\
          & =(Tu)^{-1}(y),\\
          \endalign $$
and the second boundary value problem (1.14), (1.7) is equivalent
to
$$\align
|\det\, D_y(T^*u^*)|& =g(y)/f(T^*u^*)
                 \ \ \text{in}\ \ \Om^*,\tag 2.39\\
  T^*\Om^*&=\Om.\tag 2.40\\
  \endalign $$
Noting that the defining functions $\phi$ and $\phi^*$ may be
chosen so that $\D\phi=\gamma$, $\D \phi^*=\gamma^*$ on $\pom$,
$\pom^*$ respectively, we clearly have for $x\in\pom, y\in
Tu(x)\in\pom^*$,
$$\align
\beta\cdot\gamma(x)
 &= c^{k, i}(x, y)\phi_i\phi^*_k(y)\tag 2.41\\
 &= \beta^*\cdot\gamma^*(y),\\
 \endalign$$
where
$$\beta^*(y)=D_q\phi(\A^*(D_yu^*,y)). \tag 2.42$$
Hence the quantity $\beta^*\cdot\gamma^*$ is minimized on $\pom^*$
at the point $y_0=Tu(x_0)$. Furthermore, for $y=Tu(x)$,
$x\in\pom$,
$$w^{ij}\gamma_i\gamma_j(x)
     =w^*_{kl}(y)c^{k, i}c^{l,j}(x, y)\phi_i\phi_j(x), \tag 2.43$$
where
$$w^*_{kl}(y)=u^*_{y_ky_l}(y)-c_{, kl}(x, y). \tag 2.44$$
Applying now the estimate (2.35) to $u^*$ at the point $y_0\in
\pom^*$, we finally conclude from (2.13) the desired obliqueness
estimate
$$\beta\cdot\gamma\ge \delta\tag 2.45$$
on $\pom$ for some positive constant $ \delta$ depending only on
$\Om,\Om^*, c,$ and $\B$.

The above argument clearly extends to arbitrary positive terms $B$
 (1.15). Noting also that it suffices in the above argument that
T need only be one-to-one from  a neighbourhood of the point  $x_0$
to a neighbourhood of $y_0$, we have the following theorem.

\proclaim{Theorem 2.1} Let $c\in C^3(\R^n\times\R^n)$ be a cost
function satisfying hypotheses A1, A2, with respect to bounded
$C^3$ domains $\Om, \Om^*\subset\R^n$, which are respectively
uniformly $c$-convex, $c^*$-convex with respect to each other. Let
$\B$ be a positive function in $C^1(\bom\times \R\times\R^n)$.
Then any elliptic solution $u\in C^3(\bom)$,  of the second boundary
 value problem (1.14), (1.7)
satisfies the obliqueness estimate (2.45).
\endproclaim

Note that $Tu$ is automatically globally one-to-one under
the hypotheses of Theorem 1.2 by virtue of the change of 
variables formula.  In  ensuing work , (see  [T3]), we 
 extend Theorem 2.1 to the more general
prescribed Jacobian equation (1.4). The main difference is that we
cannot directly use the $c$-transform  to get the complementary estimate
to (2.35) Instead the quantities there are transformed using the local
diffeomorphism Tu. Indeed we could also have avoided the use of duality
in the proof of Theorem 2.1 by direct transformation of (2.31). 

\vskip20pt

\centerline {\bf 3. Global second derivative bounds}

\vskip10pt

In this section we show that the second derivatives of elliptic
solutions of equation (1.14) may be estimated in terms of their
boundary values. For this estimation and the boundary estimates in
the next section, it suffices to consider the general form (1.1)
under the assumption that the matrix valued function $A\in
C^2(\Om\times \R\times\R^n)$ satisfies condition A3w, that is
$$D_{p_kp_l}A_{ij}(x, z, p)\xi_i\xi_j\eta_k\eta_l\ge 0\tag 3.1$$
for all $(x, z, p)\in \Om\times\R\times \R^n$, $\xi, \eta
\in\R^n$, $\xi\perp \eta$. We also assume $A$ is symmetric,which
is the case for the optimal transportation equation (1.14).  When (3.1) is
strengthened to the condition A3 in [MTW], that is
$$D_{p_kp_l}A_{ij}(x, z, p)\xi_i\xi_j\eta_k\eta_l
                       \ge \delta |\xi|^2|\eta|^2   \tag 3.2$$
for some constant $\delta>0$, for all $(x, z, p)\in
\Om\times\R\times \R^n$, $\xi, \eta\in\R^n$, $\xi\perp \eta$, then
the global second derivative estimate follows immediately from our
derivation of interior estimates in [MTW]. In the general case the
proof is much more complicated and we need to also assume some
kind of barrier condition, namely that there exists a function
$\wtt\phi\in C^2(\bom)$ satisfying
$$[D_{ij}\wtt\phi(x)-D_{p_k}A_{ij}(x, z, p)D_k\wtt\phi(x)]
          \xi_i\xi_j  \ge    \wtt\delta |\xi|^2\tag 3.3$$
for some positive $\wtt\delta>0$ and for all $\xi\in\R^n$,
$x,z, p\in$ some set $U\subset \Om\times\R\times\R^n$,
whose projection on $\Om$ is $\Om$. In general,
condition (3.1) implies some restriction on the domain $\Om$, but
for the case of equations arising in optimal transportation,  it can be
avoided by a duality argument. 

Our reduction to the boundary estimation follows the approach in
[GT], originating with Pogorelov, with some modification analogous
to that in [LTU]. Let $u\in C^4(\Om)$ be an elliptic solution of
equation (1.1), with $x,u(x), Du(x)\in U$ for $x\in\Om$ and $\xi$ a unit vector
in $\R^n$. Let $v$ be the auxiliary function given by
$$v=v(\cdot, \xi)=\log(w_{ij}\xi_i\xi_j)
                    +\tau|Du|^2+\kappa\wtt\phi, \tag 3.4$$
where $w_{ij}=D_{ij}u-A_{ij}$. By differentiation of equation
(1.1), we have
$$\align
 w^{ij}\big[ & D_{ij}u_\xi-D_\xi A_{ij}-(D_z A_{ij})u_\xi
                       -(D_{p_k}  A_{ij}) D_ku_\xi\big]\tag 3.5\\
 = & D_\xi\wtt B+(D_z\wtt B)u_\xi+(D_{p_k}\wtt B)D_ku_\xi, \\
 \endalign $$
where $\wtt B=\log B$. A further differentiation yields
$$\align
w^{ij} \big[D_{ij} & u_{\xi\xi}-D_{\xi\xi} A_{ij}
   -(D_{zz} A_{ij})(u_\xi)^2 -(D_{p_kp_l}A_{ij})D_ku_\xi D_lu_\xi\tag 3.6\\
  & -(D_z A_{ij})u_{\xi\xi} - (D_{p_k}A_{ij}) D_ku_{\xi\xi}
   -2(D_{\xi z}A_{ij}) u_\xi \\
  & -2(D_{\xi p_k} A_{ij}) D_ku_\xi
   -2(D_{z p_k}A_{ij})(D_ku_\xi) u_\xi \big]
   -w^{ik}w^{jl}D_\xi w_{ij}D_\xi w_{kl}\\
 = D_{\xi\xi}\wtt B & +(D_{zz}\wtt B) u_\xi^2 +
    (D_{p_kp_l}\wtt B)D_ku_\xi D_lu_\xi\\
   & +2(D_{\xi z}\wtt B)u_\xi +2(D_{\xi p_k}\wtt B) D_ku_\xi
    +2(D_{zp_k}\wtt B)(D_ku_\xi)u_\xi \\
    & +(D_z\wtt B)u_{\xi\xi} + (D_{p_k}\wtt B)D_ku_{\xi\xi.}
    \endalign $$
Furthermore differentiating (3.4) we have
$$\align
D_iv & =\frac {D_iw_{\xi\xi}}{w_{\xi\xi}}
         +2\tau D_kuD_{ki}u +\kappa D_i\wtt\phi,\tag 3.7\\
D_{ij}v&=\frac{D_{ij}w_{\xi\xi}}{w_{\xi\xi}}
      -\frac {D_iw_{\xi\xi}D_jw_{\xi\xi}}{w_{\xi\xi}^2}
      +2\tau(D_{ik}uD_{jk}u+D_k uD_{ijk}u)+\kappa D_{ij}\wtt\phi.\tag 3.8\\
\endalign            $$
where we have written $w_{\xi\xi}=D_{ij}w\xi_i\xi_j$. Using
condition A3w in (3.6) and retaining all terms involving third
derivatives, we estimate
$$\align
Lu_{\xi\xi}: & = w^{ij} (D_{ij}u_{\xi\xi}+b^{ij}_kD_ku_{\xi\xi})
        -(D_{p_k}\wtt B)D_ku_{\xi\xi}\tag 3.9\\
    &\ge w^{ik}w^{jl}D_\xi w_{ij}D_\xi w_{kl}
        -C\{(1+w_{ii})w^{ii}+(w_{ii})^2\}\\
        \endalign $$
where, as in the previous section, $b^{ij}_k=-D_{p_k}A_{ij}$ and
$C$ is a constant depending on the first and second derivatives of
$A$ and $\log B$ and $\sup_\Om (|u|+|Du|)$. To apply A3w, we fix a
point $x\in\Om$ and choose coordinate vectors as the
eigenfunctions of the matrix $[w_{ij}]$ corresponding to
eigenvalues $0<\lam_1\le \cdots\le \lam_n$. Writing $A_{ij,
kl}=D_{p_kp_l}A_{ij}$, we then estimate
$$\align
w^{ij}A_{ij, kl} u_{k\xi}u_{l\xi}
  &\ge w^{ij}A_{ij, kl} w_{k\xi}w_{l\xi} -Cw^{ii}(1+w_{ii})\\
  &\ge \sum_{k\,\text{or}\, l=r} \frac {1}{\lam_r} A_{rr, kl}
               (\lam_k\xi_k)(\lam_l\xi_l)-Cw^{ii}(1+w_{ii})\\
  &\ge -C\{w^{ii}(1+w_{ii})+w_{ii}\}\\
  \endalign $$
From (3.9), we obtain also
$$
Lw_{\xi\xi}\ge w^{ik}w^{jl}D_\xi w_{ij}D_\xi w_{kl}
        -C\{(1+w_{ii})w^{ii}+ w_{ii}^2\} \tag 3.10$$
for a further constant $C$. Here we use equation (3.5) to control
the third derivative term arising from differentiating
$A_{kl}\xi_k\xi_l$. From (3.8) and (3.10), we obtain, after some
reduction,
$$\align
Lv\ge & \frac{1}{w_{\xi\xi}} w^{ik}w^{jl}D_\xi w_{ij}D_\xi w_{kl}
   -\frac {1}{w_{\xi\xi}^2} w^{ij}D_i w_{\xi\xi}D_jw_{\xi\xi}\tag 3.11\\
   &+2\tau w_{ii} +\kappa w^{ii} -C\{ \frac
   {1}{w_{\xi\xi}}[(1+w_{ii})w^{ii}+w_{ii}^2]+\tau+\kappa\}\\
   \endalign $$
Now suppose $v$ takes its maximum at a point $x_0\in\Om$ and a
vector $\xi$, which we take to be $e_1$. We need to control the
first two terms on the right hand side of (3.11). To do this we
choose remaining coordinates so that $[w_{ij}]$ is diagonal at
$x_0$. Then we estimate
$$\align
  &\frac{1}{w_{\xi\xi}} w^{ik}w^{jl}D_\xi w_{ij}D_\xi w_{kl}
   -\frac {1}{w_{\xi\xi}^2} w^{ij}D_i w_{\xi\xi}D_jw_{\xi\xi}\tag 3.12\\
 &= \frac{1}{w_{11}} w^{ii}w^{jj}(D_1 w_{ij})^2
   -\frac {1}{w_{11}^2} w^{ii}(D_i w_{11})^2\\
 &\ge  \frac {1}{w_{11}^2}
 \sum_{i>1}[2w^{ii}(D_1w_{1i}^2-w^{ii}(D_iw_{11})^2]\\
&= \frac {1}{w_{11}^2} \sum_{i>1}w^{ii}(D_iw_{11})^2 \\
   & \qquad  +\frac {2}{w_{11}^2} \sum_{i>1}w^{ii}
    [D_1w_{1i}-D_iw_{11}]\, [D_1w_{1i}+D_iw_{11}]\\
  & \ge \frac {1}{w_{11}^2} \sum_{i>1}w^{ii}(D_iw_{11})^2 \\
   & \qquad +\frac {2}{w_{11}^2} \sum_{i>1}w^{ii}
    [D_iA_{11}-D_1A_{1i}]\, [2D_iw_{11}+D_iA_{11}-D_1A_{1i}]\\
  & \ge -Cw^{ii}\\
 \endalign $$
Combining (3.11) with (3.12), we obtain the estimate, at $x_0$,
$$ Lv\ge \tau w_{ii} +\kappa w^{ii} - C\{\tau+\kappa\}, \tag 3.13 $$
for either $\tau$ or $\kappa$ sufficiently large. Note that when we use
(3.7) in the second last line of (3.12), we improve (3.13) by retention of the term
$\sum_{i>1}w^{ii}(D_iw_{11})^2$ on the right hand side , which 
corresponds to the key term in the Pogorelov argument for interior
 estimates [GT].

From (3.13), we finally obtain an estimate from above for $w_{ii}(x_0)$,
which we formulate in the following theorem.

\proclaim {Theorem 3.1} Let $u\in C^4(\Om)$ be an elliptic
solution of equation (1.1) in $\Om$ , with  $x,u(x),Du(x) \in U$,
for all $x\in\Om$. Suppose
the conditions A3w and (3.3) hold and $B$ is a positive function
in $C^2(\bom\times\R\times\R^n)$. Then we have the estimate
$$\sup_\Om |D^2 u|\le C(1+\sup_{\pom} |D^2 u|), \tag 3.14$$
where the constant $C$ depends on $A, B, \Om$,$\Om^*$ and
$\sup_U(|z|+|p|)$.
\endproclaim

Note that we only need the condition A3w to hold on the set $U$.

From the proof of Theorem 3.1 we obtain the corresponding estimate
for equation (1.14), without the barrier condition (3.3).

\proclaim {Theorem 3.2} Let $u\in C^4(\Om)$ be an elliptic
solution of equation (1.14) in $\Om$ with  $Tu(\Om)\subset\Om^*$.
Suppose the cost function $c$ satisfies hypotheses A1, A2, A3w and
$B$ is a positive function in $C^2(\bom\times\R\times\R^n)$. Then
we have the estimate (3.14).
\endproclaim

To prove Theorem 3.2, we take $\kappa=0$ in the proof of Theorem
3.1, to obtain an estimate for  $w_{ii}$ in terms of $w^{ii}$,
that is
$$w_{ii}\le \eps\sup_\Om w^{ii}+C_\eps \big ( 1 + \sup_{\partial \Omega}
 |D^2u| \big ), \tag 3.15$$
for arbitrary $\eps>0$, with constant $C_\eps$ also depending on
$\eps$. If $T$ is globally one-to-one, we then conclude (3.13), in
the optimal transportation case, by using the dual problem (2.37),
(2.38). More generally, we consider the dual
function $v^*$ in place of (3.4), given by

$$
 v^* \, = \, v(x,\xi)
 \, = \, \log \big ( w^{ij}
  c_{i,k} c_{j,l}
  \xi_k \xi_l
  \big )
   \, + \, \tau |c_y(x,Tu(x))|^2 \tag 3.16
$$

\noindent and suppose it is maximized at a point $x^*_0$ in
$\Omega$. Since $T$ will now be one-to-one from a neighborhood $\Cal N$
of $x^*_0$ to a neighbourhood $\Cal N^*$ of $y^*_0 = Tu(x^*_0)$, we may then
 proceed as before, noting that in $\Cal N^*$, $v^*$ is given by (3.4) with $u$
replaced by its $c-$transform $u^*$.

The estimate (3.14) arose from our investigation of the classical
Dirichlet problem , (see [T3] ). We remark also that from (3.16), we
see that  (3.3) is also not needed when $n=2$.

 \vskip20pt

\centerline{\bf 4. Boundary estimates for the second derivatives}

\vskip10pt

This part of our argument is similar to the treatment of the
oblique boundary value problems for Monge-Amp\`ere equations in
[LTU, U3]. The paper [LTU] concerned the Neumann problem,
utilizing a delicate argument which did not extend to other linear
oblique boundary conditions. For nonlinear oblique conditions of
the form (2.1) where the function $G$ is uniformly convex in the
gradient, the twice tangential differentiation of (2.1) yields
quadratic terms in second derivatives which compensate for the
deviation of $\beta=G_p$ from the geometric normal and permit some
technical simplification for general inhomogeneous terms $\B$
[U3].

First we deal with the non-tangential second derivatives. Letting
$F\in C^2(\bom\times\R\times\R^n)$ and $v=F(\cdot, u, Du)$, where
$u\in C^3(\bom)$ is an elliptic solution of equation (1.1), we
have from our calculation in Section 2,
$$|Lv|\le C(w^{ii}+w_{ii}+1),\tag 4.1$$
where $L$ is given by (2.17) and $C$ is a constant depending on
$A, B, G, \Om$ and $|u|_{1; \Om}$. Now using the equation (1.1)
itself, we may estimate
$$w_{ii}^{\frac{1}{n-1}} \le Cw^{ii}, \tag 4.2$$
so that, writing $M=\sup_\Om w_{ii}$, we have from (4.1)
$$|Lv|\le C(1+M)^{\frac {n-2}{n-1}} w^{ii}. \tag 4.3$$
Hence, if there exists a $C^2$ defining function $\phi$ satisfying
(3.3) near $\pom$, together with $\phi=0$ on $\pom$, we obtain by
the usual barrier argument, taking $F=G$,
$$|D(\beta\cdot Du)|\le C(1+M)^{\frac{n-2}{n-1}}\tag 4.4$$
on $\pom$, so that in particular
$$w_{\beta\beta} \le C(1+M)^{\frac {n-2}{n-1}}\tag 4.5$$
on $\pom$. Now for any vector $\xi\in\R^n$, we have
$$w_{\xi\xi}=w_{\tau\tau}
    +b(w_{\tau\beta}+w_{\beta\tau})+b^2w_{\beta\beta},\tag 4.6$$
where
$$ b=\frac {\xi\cdot\gamma}{\beta\cdot \gamma},
 \ \ \ \tau=\xi-b\beta. \tag 4.7$$
Suppose $w_{\xi\xi}$ takes its maximum over $\pom$ and tangential
$\xi$, $|\xi|=1$ at $x_0\in\pom$ and $\xi=e_1$. Then from (4.5)
and (4.6) and tangential differentiation of the boundary condition
(2.1) we have on $\pom$,
$$w_{11}\le |e_1-b\beta|^2 w_{11}(0)+bF(\cdot, u, Du)
   +Cb^2(1+M)^{\frac{n-2}{n-1}},\tag 4.8$$
for a given function $F\in C^2(\bom\times\R\times\R^n)$. Combining
(2.26), (3.10), (4.1) and (4.2) and utilizing a similar barrier argument
to that is Section 2, we thus obtain the third
derivative estimate
$$-D_\beta w_{11}(x_0)\le C(1+M)^{\frac {2n-3}{n-1}}, \tag 4.9$$
Differentiating (2.1) twice in a tangential direction $\tau$, with
$\tau(x_0)=e_1$, we obtain at $x_0$,
$$(D_{p_kp_l} G)u_{1k}u_{1l}+
  (D_{p_k} G)u_{11k}\le C(1+M), \tag 4.10$$
  whence we conclude from (4.9)
$$\max_{\pom} |D^2 u|\le
        C(1+\sup_\Om |D^2 u|)^{\frac {2n-3}{2(n-1)} }\tag 4.11$$
by virtue of the uniform convexity of $G$ with respect to $p$.
Taking account of the global estimate (3.13), we complete the
proof of Theorem 1.1. $\square$

\vskip10pt

Once the second derivatives are bounded, the equation (1.1) is
effectively uniformly elliptic so that from the obliqueness
estimate (2.43), we obtain global $C^{2,\alpha}$ estimates from
the theory of oblique boundary value problems for uniformly
elliptic equations in [LT]. By the theory of linear elliptic
equations with oblique boundary conditions [GT], we then infer
estimates in $C^{3, \alpha}(\bom)$ for any $\alpha<1$ from the
assumed smoothness of our data. We may also have assumed that our
solution $u\in C^2(\bom)$.

As in the previous section, the technicalities are simpler when
condition A3w is strengthened to condition A3 and we also obtain
local boundary estimates for the second derivatives. To see this
we estimate the tangential second derivatives first by
differentiating the equation (1.1) and boundary condition (2.1)
twice with respect to a tangential vector field $\tau$ near a
point $y\in \pom$. We then obtain an estimate for $\eta
D_{\tau\tau} u$, for an appropriately chosen cut-off function
$\eta$. The mixed tangential-normal second derivatives $D_{\tau
n}u$ are estimated as above by a single tangential differentiation
of (2.1) so that the double normal derivative may be obtained
either from (4.5) or from the equation (1.1) itself and the
estimates in Section 2 for $w^{ij}\gamma_i\gamma_j$ from below,
similarly to the Dirichlet problem, see [T1].

\vskip20pt

\centerline{\bf 5. Method of continuity}

\vskip10pt

To complete the proof of Theorem 1.2, we adapt the method of
continuity for nonlinear oblique boundary value problems, presented
in [GT] and already used in the special case (1.16) (1.17) [U1].
The situation here is  more complicated 
 unless we know in advance that there exists a smooth
function $u_0$, satisfying the ellipticity condition (1.2) together
with the boundary condition (1.7).Later in Section 7,
we shall prove the existence of such a function, (at least
for approximating domains). Otherwise we need to consider
families of subdomains. To commence the procedure,
 we fix a point $x_0\in\Om$. Then for sufficiently small radius
$r>0$, the ball $\Om_0=B_r(x_0)\subset \Om$ will be uniformly
$c$-convex with respect to $\Om^*$ and the function $u_0$, given
by
$$u_0(x)=\frac \kappa 2 |x-x_0|^2+p_0\cdot (x-x_0), \tag 5.1$$
will satisfy the ellipticity condition (1.2). Moreover the image
$\Om_0^*=Tu_0(\Om_0)$ will be uniformly $c^*$-convex with respect
to $\Om$ with $Tu_0$  a diffeomorphism from $\Om_0$ to
$\Om_0^*$. To see this we observe that
$$c_x(x_0, \Om_0^*) =B_{\kappa r}(p_0), \tag 5.2$$
so that by taking $\kappa r$ small enough, we can fulfill
condition (1.19) on $\pom_0^*$, with respect to $x_0\in\Om$, for
constant $\delta_0=\frac {1}{\kappa r}$ as large as we wish. Suppose
now we can foliate $\Om -\Om_0$ and $\Om^* -\Om_0^*$ by boundaries
of  $c$-convex and $c^*$- convex domains, respectively. That is there
exist increasing families of domains $\{\Om_t\}$,
$\{\Om^*_t\}$, $0\le t\le 1$, continuously depending on the parameter $t$,
 such that
\newline
(i) $\Om_t\subset\Om$, $\Om_t^*\subset\Om^*$,
\newline
(ii) $\Om_1=\Om$, $\Om_1^*=\Om^*$,
\newline
(iii) $\pom_t, \pom^*_t\in C^4$,  uniformly with respect to $t$,
\newline
(iv) $\Om_t, \Om^*_t$ are uniformly $c$-convex, $c^*$-convex with
respect to $\Om^*, \Om$, respectively.
\newline
The construction of such a family is discussed below.

Given our families of domains $\Om_t, \Om^*_t$, $0<t\le
1$, we need to define corresponding equations. Let $B$ be a
positive function in $C^2(\bom\times \R^n)$ and $f$ a positive
function in $C^2(\bom)$ such that
$$f=-\sigma u_0+\log [ \det \{D^2 u_0-c_{xx}(\cdot, \A(\cdot,Du_0)\}/
 B(\cdot, Du_0)]\tag 5.3$$
in $\Om_0$, for some fixed constant $\sigma>0$. We then consider
the family of boundary value problems:
$$\align
F[u]:&=\det \{D^2 u-c_{xx} (\cdot, \A(\cdot, Du)\}
    = e^{\sigma u+(1-t)f} B (\cdot, Du),\tag 5.4\\
Tu(\Om_t) & = \A(\cdot, Du)(\Om_t)=\Om_t^*, \\
\endalign $$
From our construction and the obliqueness, we see that $u_0$ is
the unique elliptic solution of (5.4) at $t=0$.

From Section 2, we also see that the boundary condition in (5.4)
is equivalent to the oblique condition
$$G_t(\cdot, Du):= \phi_t^* (Y(\cdot, Du)) =0
                  \ \ \ \text{on}\ \ \pom_t. \tag 5.5$$
To adapt the method of continuity from [GT], we fix $\alpha\in (0,
1)$ and let $\Sigma$ denote the subset of $[0, 1]$ for which the
problem (5.4) is solvable for an elliptic solution $u=u_t\in C^{2,
\alpha}(\bom_t)$, with $Tu$ invertible. We then need to show that
$\Sigma$ is both closed and open in $[0, 1]$. First we note that
the boundary condition (5.4) implies a uniform bound for $Du_t$.
Integrating the equation (5.4), we then obtain uniform bounds for
the quantities
$$\int_{\Om_t} e^{\sigma u_t}, $$
so that the solutions $u_t$ will be uniformly bounded for
$\sigma>0$. Uniform estimates in $C^{2,1}(\bom)$ now follow from
our a priori estimates in Section 4, which are also clearly
independent of $t\in [0,1]$. By compactness, we then infer that
$\Sigma$ is closed. To show $\Sigma$ is open, we use the implicit
function theorem and the linear theory of oblique boundary value
problems, as in [GT]. The varying domains $\{\Om_t\}$ may be
handled by means of diffeomorphisms approximating the identity,
which transfer the problem (5.4) for $t$ close to some
$t_0\in\Sigma$ to a problem in $\Om_{t_0}$. We then conclude the
solvability of (5.4) for all $t\in [0, 1]$, which implies there
exists a unique elliptic solution $u=u_\sigma\in C^3(\bom)$ of the
boundary value problem
$$\align
F[u] &= e^{\sigma u} B(\cdot, Du), \ \  \ \tag 5.6\\
Tu(\Om) & =\Om^*\\
\endalign $$
for arbitrary $\sigma>0$, with $Tu$ one-to-one. To complete the
proof of Theorem 1.2,(at least when the above foliations exist),
 we assume that $B$ satisfies (1.5), (1.8)
and (1.9). As above we see that the integrals
$$\int_\Om e^{\sigma u_\sigma}$$
are uniformly bounded, with $D(\sigma u_\sigma)\to 0$ as
$\sigma\to 0$. Consequently $\sigma u_\sigma\to\,$constant$\,=0$
by (1.19) and modulo additive constants, $u_\sigma\to u$ as
$\sigma\to 0$, where $u$ is the solution of (1.14), (1.7), as
required. 

We may construct the family of domains  $\{\Om_t\}$ used above, 
if we are given a defining
function $C^4$ defining function $\phi$, satisfying

$$[D_{ij}\phi(x)-c^{l,k}c_{ij, l}(x,y)D_k\phi(x)]\xi_i\xi_j
  \ge \delta_0 |\xi|^2\tag 5.7$$
for all $x\in\Om, y\in \Om^*$, $\xi\in \R^n$, which takes its minimum
  at  $x_0$. Note that the uniform $c$-convexity of $\Om$ implies
  the existence of a defining function satisfying (5.7) in
  a neighbourhood of  $\pom$, as in (2.26). There are various ways of
  constructing suitable families from a global defining function, $\phi$.
  In particular taking $\phi (x_0) = -1$, we may choose 

  $$\Om_t=\{x\in\Om\ |\ \phi_t(x)<0\}$$
where  $\phi_t$ is defined by

  $$\phi_t=t (\phi -a)+(t_0-t)\phi_0, \tag 5.8$$
for $\phi_0(x)=|x-x_0|^2-r^2$ , $t \le t_0$, for some $0<t_0<1$ 
 and $a$ close to $-1$ to first deform to a small sub-level set of $\phi$,
 followed by taking $\phi _t = (1-t)a/(1-t_0)$ for $t \ge t_0$. Alternatively,
 we could have chosen $\Om_0 = { \phi<a}$ at the outset and only used
the second deformation.The domains $ \Om^*_t$
  may be similarly constructed. If the curvatures of   $\pom^*$
 are sufficiently large, for example if $\Om^*$ is a small ball, then the
 existence of a defining function satisfying (5.7) follows by pulling back
 from a single image $c_y(\cdot,y_0)(\Om)$. As a byproduct, we see that
 if $\Om$ is uniformly $c$-convex with respect to a single point
 $y_0\in\Om^*$, there at least exists a smooth function $u_0\in C^3(\bom)$
 satisfying the ellipticity condition (1.2). Moreover $Tu_0$ is a diffeomorphism 
 from $\Om$ to  a small ball $B_r(y_0)\subset\Om^* $.
 
 From the above  considerations, we see that the proof of Theorem 1.2
 is completed in the cases where either $\Om$ or $\Om^*$ is a small ball.
 The general case will then follow by further use of the method of continuity
 if  there exists a defining function satisfying (5.7) for either domain. However
 we will take up a different approach from this point and use the function $u_0$
 constructed above to construct a further function $u_1$ approximately
 satisfying our given boundary conditions to which the method of continuity can
 be applied without domain variation. Specifically we will prove
 
 \proclaim{Lemma 5.1} Let the domains $\Om$ and $\Om^*$ and cost function 
 $c$ satisfy the hypotheses of Theorem 1.1. Then for any $\epsilon > 0$, there
  exists a uniformly   $c^*$-convex approximating domain, $\Om^*_\epsilon$,
   lying within distance $\epsilon$ of  $\Om^*$,  and satisfying the corresponding
condition (1.19) for fixed $\delta_0$, together with a function
 $u_1\in C^4(\bom)$ satisfying the ellipticity condition (1.2) and the 
 boundary condition (1.7) for  $\Om^*_\epsilon$ .
\endproclaim

From Lemma 5.1 we complete the proof of Theorem 1.2. We defer the proof of
 Lemma 5.1 to Section 7 as the proof will use some of the same 
ingredients as in our discussion of optimal transportation in the next section.
In Section 7,  we will also indicate an alternative  and direct
construction of the function
$u_1$, which avoids domain variation altogether in the method of continuity.

\vskip20pt

\vskip20pt

\vskip20pt

\centerline {\bf  6. Optimal Transportation}

\vskip10pt

The interior regularity of solutions to the optimal transportation
problem is considered in [MTW],
 under conditions A1, A2, A3 and the $c^*$- convexity of the
target domain  $\Omega^*$.  Our approach is to first show that the
Kantorovich  potentials are generalized solutions of the boundary
value  problem (1.14), (1.7) in the  sense of Aleksandrov and
Bakel$'$man. The $c^*$-convexity of $\Omega^*$ is used to show
the image of the generalized normal mapping lies in
$\overline{\Omega}^*$ and condition A3 is employed to obtain a
priori  second derivative estimates from which the desired
regularity  follows. The potential functions $u$ and $v$ solve the
dual  problem of minimizing the functional
$$I(u,v)= \int _{\Omega} f u  +  \int _{\Omega^*} g  v \tag 6.1$$
over the set $K$ given by
$$K = \big \{ (u,v) \ \big | \ u, v\in C^0(\Omega), C^0(\Omega^*)
  \ \ \text{resp.}  \ \ u(x)+v(y) \ge c(x,y)
 \ \text{for all}\ x \in \Omega , y \in \Omega ^* \big \}. \tag 6.2$$
The potential functions $(u,v)$ satisfy the relations
$$\align
 u(x) &= \sup_{y\in \Omega}\big \{ c(x,y)-v(y) \big \}, \tag 6.3\\
 v(y) &= \sup_{x\in \Omega} \big \{ c(x,y)-u(y) \big \},\\
\endalign $$
\noindent that is they are the $c^*$ and $c$ transforms of each
other. Since $c \in C^{1,1},$ they will be semi-convex. The
optimal mapping $T$ is then given almost everywhere by (1.23) and
the equation (1.14) will be satisfied with elliptic solution $u$
 almost everywhere in $\Omega$. The functions $u$ and $v$
are respectively $c$ and $c^*$- convex. A function $u \in
C^0(\Omega)$ is called $c-convex$ in $\Omega$ if for each $x_0 \in
\Omega$, there exists $y_0 \in \Bbb{R}^n$ such that
$$ u(x) \ge u(x_0) + c(x, y_0) -c(x_0, y_0) \tag 6.4$$
for all $x \in \Omega$. If $u$ is a $c-$convex function, for which
the mapping $T$ given by (1.23) is measure preserving, then it
follows that $u$ is a potential and again $T$ is the unique
optimal mapping. These results hold under the hypotheses A1 and
A2 and it suffices to assume the densities $f,g \ge 0, \in L^1(\Om),
L^1(\Om^*)$, respectively, whence the mapping T is only determined
almost everywhere on the set where f is positive. The reader is 
referred to [C3,GM,MTW,U2,V] for further details.

From the above discussion we see that the solution of the boundary
value problem (1.14), (1.7) will automatically  furnish a
potential for the optimal transportation problem if it is
$c$-convex. Note that ellipticity only implies that the solution
is locally $c$-convex and we need a further argument to conclude
the global property, unlike the case of quadratic cost functions
and convex solutions. First we recall the concept of generalized
solution introduced in [MTW]. Let $u$ be a $c$-convex function on
the domain $\Omega$. The $c-normal$ mapping, $\chi_u$, is
defined by
$$ \chi_u(x_0) = \big \{ y_0 \in \Bbb{R}^n \ \big | \
 u(x) \ge u(x_0)+c(x,y)-c(x_0,y_0) , \text{ for all } x \in \Omega
 \big \}. \tag 6.5$$
Clearly, $\chi_u (x_0) \subset Y(x_0, \partial u(x_0) )$ where
$\partial$ denotes the subgradient of $u$. For $g \ge 0, \in
L^1_{\ell {\text oc}}(\Bbb{R}^n)$, the generalized
\it{Monge-Amp\`ere measure} $\mu [u,g]$ \rm{is then defined by}
$$ \mu [u,g] (e) =  \int_{\chi_u (e)} g \tag 6.6 $$
for any Borel set $e \in \Omega$, so that $u$ satisfies equation
(1.14) in the generalized sense if
$$ \mu [u,g] = f \, \text{d} x \, . \tag 6.7 $$
The boundary condition (1.7) is satisfied in  the generalized
sense if
$$ \Omega^* \subset \chi_u(\overline{ \Omega}), \qquad
 \big | \big \{ \, x \in
  \Omega \ \big | \ f(x) > 0 \text{ and } \chi _u (x)  -
  \overline{\Omega}^* \ne \emptyset \big \} \big | = 0 \tag 6.8 $$
The theory of generalized solutions replicates that for  the
convex case, $c(x,y) = x \cdot y$, [MTW]. If $f$ and $g$ are
positive, bounded measurable functions on $\Omega, \Omega^*$
respectively satisfying the mass balance condition (1.9), and  $c$
satisfies A1, A2, then there exists a unique (up to  constants)
generalized solution of (6.7), (6.8), (with $g=0$  outside
$\Omega^*$), which together with its $c$ transform $v$,  given by
(6.3), uniquely solves the dual problem (6.1),  (6.2)([MTW]).

Now let $u\in C^2(\overline{\Omega})$ be an elliptic solution of
the boundary value problem (1.7), (1.14) and $v$ a $c-$convex
solution of the corresponding generalized problem. By adding
constants, we may assume $\inf_{\Omega}(u-v)=0$. We need to prove
$u=v$ in $\Omega$, that is the strong comparison principle holds.
Let $\Omega'$ denote the subset of $\Omega$ where $u>v$ and first
suppose that $\partial   \Omega' \cap\Omega \ne \emptyset$. Note
that if  $v\in C^2(\Omega)$, this
situation is immediately ruled out by the classical strong maximum
principle [GT]. Otherwise we may follow the proof of the strong
maximum principle as there will exist a point $x_0 \in \partial
\Omega' \cap \Omega$, where $\Omega'$ satisfies an interior sphere
condition, that is there exists a ball $B \subset \Omega -
\Omega'$ such that $x_0 \in\partial \Omega' \cap
\partial B$, $u(x_0) = v(x_0)$ and $u>v $ in $B$. Since $v$
is semi-convex, $v$ will be twice differentiable at $x_0$, with
$Dv(x_0) = Du(x_0)$. Moreover by passing to a smaller ball if
necessary we may assume both $u$ and $v$ are $c-$convex in $B$.
Since $u$ is a  smooth elliptic solution of (1.14), there will
exist a strict supersolution $w \in C^2(\overline{B}-B_{\rho})$,
for some concentric ball $B_{\rho}$ of radius $\rho < R$,
satisfying $w(x_0) = u(x_0)$, $w \ge v$ on $\partial B \cup
\partial B_\rho, Dw(x_0) \ne Du(x_0)$. By the comparison
principle, [MTW], Lemma 5.2, we have $w\ge v$ in $B-B_{\rho}$, and
hence $Dw(x_0) = Dv(x_0)$, which is a contradiction. Thus we may
assume $\partial \Omega' \cap \Omega = \emptyset$, that is $u>v$
in $\Omega$ with $u(x_0) = v(x_0)$ for some point $x_0 \in
\partial \Omega$. From our argument above, we obtain a
function $w \in C^2(\overline{B} -B_{\rho})$ satisfying $w(x_0) =
u(x_0) = v(x_0)$, $v\le w \le u$ in  $ \overline{B} - B_{\rho}$,
together with
$$ u(x) - w(x)  \ge \epsilon |x-x_0| \tag 6.9$$
for all $x \in B_R - B_{\rho}$. Since $v \le w$ in $B_R -
B_{\rho}$ , this contradicts  the
obliqueness  condition (2.43) if $\Omega^*$ is $c^*-$convex.

Alternatively we may proceed directly ( and more simply) as follows
 to show that the
solution $u$ is c-convex, using the property that $Tu$ is one-to-one.
In fact, as mentioned previously in Section 2, this would follow 
automatically from the change of variables 
formula by virtue of  the mass balance condition (1.9).
 Let $x_0\in\Omega$ and $y_0=Tu(x_0)$. Suppose there exists a point
$x_1\in\Omega$, where
$$ u(x_1) < c(x_1, y_0) - c(x_0, y_0).\tag 6.10$$
By downwards vertical translation, there exists a point $x_2 \in
\partial\Omega$, satisfying
$$ u(x) > u(x_2) + c(x,y_0) - c(x_2, y_0).\tag 6.11$$
for all $x \in \Omega$. Putting $y_2=Tu(x_2)$, we must also have
$$c_x(x_2,y_2).\gamma(x_2) < c_x(x_2,y_0).\gamma(x_2),\tag6.12$$
which again contradicts the $c^*-$convexity of $\Omega^*$.

\vskip10pt

{\bf Remarks.}

\vskip10pt

In the first proof above, we employed a comparison result that if $u$ is a
 classical elliptic supersolution of (1.14) dominating a generalized
 subsolution $v$ on the boundary of a subdomain $\Om' $, then $u\ge v$
 in $\Om'$.  In our local uniqueness argument in [MTW], we also used implicitly
 the complementing result that if u is an elliptic subsolution dominated by a
 generalized supersolution $v $ on $\partial \Omega'$, then $u\le v$ in $\Om'$.
 However,in this case,we cannot apply Lemma 5.2 in [MTW] directly as
 local $c$-convexity  of  $v$ may not imply global $c$-convexity in $\Om$, unless
 $v\in C^1(\Om)$. This situation is rectified in [TW], under the strong A3 hypothesis.
 However if $\Om$ and $\Om^*$ lie respectively
 in domains $\Om_0$ and $\Om^*_0$ satisfying the hypotheses of 
 Theorem1.2, with  $\Om^*$ also $c^*$-convex with respect to $\Om_0$,
 and f and g are positive in $L^1(\Om)$ and $L^1(\Om^*)$ respectively,
 then it follows directly by approximation from Theorem 1.3 that the local 
 $c$-convexity of the potential u solving the Kantorovich dual problem
 implies its global $c$-convexity.   Other results, such as the $c$-convexity 
 of the contact set under condition
 A3w, also follow from Theorem 1.3 by approximation. The reader is
 referred to Loeper [L]  for a full treatment of this approach, including
 the sharpness of condition A3w for regularity.

\vskip20pt

\centerline {\bf  7. Completion of Proof of Theorem 1.2}

\bigskip

In this section, we provide the proof of Lemma 5.1, thereby completing
that of Theorem 1.2. For this purpose, we need to draw on a geometric 
property of $c$-convex domains introduced in [TW]. Namely, suppose
that  $\Om$ is uniformly $c$-convex, with respect to $\Om^*$, and 
that the cost function $c$ satisfies
conditions A1, A2 and A3w. Denoting as before the unit outer normal to
$\partial \Omega$ by $\gamma$, we see that the level set $\Cal E$ of the
function $e$, given by
$$e(x) = e_y(x) = c(x,y)  - c(x,y_0)  , \tag 7.1$$
passing through $x_0$, is tangential to $\partial \Omega$ at $x_0$ if 
$$y = Y(x_0,p_0 + t\gamma_0), \tag 7.2$$
for $t > 0$, $p_0 = c_x(x_0,y_0)$, $\gamma_0 = \gamma(x_0)$
 that is, $y$ lies on the $c^*$-segment 
which is the image under $Y(x_0,\cdot)$ of the line from $p_0$ with 
slope $\gamma _0$. Then it follows from [TW] that $\Om$ lies strictly
on one side of $\Cal E$, whence
$$c(x,y) -c(x_0,y) < c(x,y_0) - c(x_0,y_0),\ \ \ x\in \Om . \tag 7.3$$
To prove (7.3) directly from (1.19), we take $x_0 = 0$, 
set $x^\prime =(x_1, .....x_{n-1})$ and choose coordinates so
that $\gamma_0 = (0, ......, -1)$. By Taylor's formula,
$$e(x) - e(x_0) \le  -tx_n + [A_{ij}(0,p_0 + t\gamma_0) -
A_{ij}(0,p_0)]x^\prime {_i}x^\prime{_j} + tO(|x||x_n| + |x|^3).$$
Using (1.13), (1.19), condition A3w and again Taylor's formula, we have
$$\align
[A_{ij}(0,p_0 + t\gamma_0) - A_{ij}(0,p_0)]x^\prime {_i}x^\prime{_j}
& \le -tc^{l, n} c_{ij, l}(x_0, y_0)x^\prime {_i}x^\prime{_j} \\
& \le tD_i\gamma_j(x_0)x^\prime {_i}x^\prime{_j} -t \delta_0|x^\prime|^2,\\
\endalign$$
so that, 
$$e(x) - e(x_0)  \le  -tx_n + tD_i\gamma_j(x_0)x^\prime {_i}x^\prime{_j}
+ tO(|x||x_n| + |x|^3)  < 0,$$
for $x \in \bom - {x_0}$, sufficiently small. Consequently $\Om$ lies
locally, strictly on one side of $\Cal E$. We can then verify 
the global inequality
(7.3) by contradiction, as in [TW]. For if (7.3) is violated, the set
$$U_a = \{ x\in \pom  ; e(x) > a|\}$$
contains two disjoint components, for sufficiently small $a>0$.
Increasing $a$, we see that these components will meet first at a point
$x^*\in \pom$ at which the level set of the function $e$ is tangential,
contradicting the local inequality at $x^*$.

Now, to commence the proof of Lemma 5.1, we take $u_0$ to be
 a function as constructed in Section 6, that is 
$u_0$ is a smooth uniformly $c$-convex  function on $\Omega$, whose
$c$-normal mapping $Tu_0 = Y(\cdot,Du_0)$ has image $\omega^*$, 
which is a $c^*$-convex
subdomain of $\Omega^*$. Here we call a $c$-convex function
uniformly c-convex if it also satisfies the ellipticity condition (1.3).
We remark that the $c$-convexity of $u_0$ could also have been proved
from (7.3), using the $c$-convexity of $\Om$ and condition A3w,
instead of the $c^*$-convexity of $\Om^*$  which we used in Section 6.
Also by approximation, we may assume $u_0\in C^\infty(\bom)$.
 A function $h$ is called a $c$-function if
it has the form
$$h(\cdot)= c(\cdot, y_0)+a$$
for $y\in\R^n$ and some constant $a$. When $c(x, y)=x\cdot y$, a
$c$-function is a linear function. Obviously the $c$-normal mapping of
$h$ is  the constant map, $Th(x)=y_0$ for all $x\in\R^n$. Let
$$u_1(x) = \sup\{u_0, h(x)\},\ \ \ x\in \Om^\delta,\tag 7.4$$
where  the sup is taken over the set $\Cal S$ of 
c-functions $h$ with $h\le u_0$ in $\Om$,
$T_h(\Om)\subset\Om^*$ and $\Om^\delta=\{x\in\R^n:\ \ \dist(x, \Om)<\delta\}$,
for some $\delta > 0$, is a neighbourhood of $\Om$. The following
lemma describes the properties of $u_1$.

\vskip10pt

\proclaim{ Lemma 7.1}. Assume that the cost function $c$ satisfies A1, A2, A3w
and the domains $\Om$ and $\Om^*$ are uniformly $c$-convex with
respect to each other. Then, for sufficiently small $\delta$,
 the function $u_1$ is a $c$-convex extension
of $u_0$ to  $\Om^\delta$,
whose $c$-normal image under $u_1$ is $\bom^*$. Moreover for any point
$x\in \Om^\delta -\bom$, there exist unique points, $x_b\in \pom$, 
$y_b\in \pom^*$, such that  $\chi_{u_1} = y_b$ on the $c$-segment, $\ell_{x_b}$,
joining $x$ to $x_b$, with respect to $y_b$,(except at the endpoint $x_b$)
with the resultant mappings 
being $C^2$ diffeomorphisms from $\pom^r$ to  $\pom,  \pom^*$ respectively,
for any $r < \delta$.
\endproclaim

\vskip10pt

Lemma 5.1 will follow from Lemma 7.1 by modification of $u_1$ outside
$\Om$  and mollification. To prove Lemma 7.1, we first take any c-function
in the set  $\Cal S$, with $c$-normal image $y\in \Om^* -\omega^*$ and increase
it until its graph meets that of $u_0$ on $\bom$ at a point $x_b$, which
will lie in $\pom$, since $u_0$ is uniformly $c$-convex. Accordingly we
obtain a $c$-function $h\in \Cal S$, given by
$$ h(x) = h_{x_b,y}(x) = c(x, y) - c(x_b, y) + u_0(x_b).\tag 7.5$$
Since $h\le u_0$ in $\Om$ and $h(x_b) = u_0(x_b)$, we see that
the point $y$ must lie on the $c^*$-segment  $\ell^*_{x_b}$,with respect
to $x_b$, starting at $y_{0,b} = Tu_0(x_b)$ and given by
$$c_x(x_b, \ell^*_{x_b}) = \{Du_0(x_b) + t\gamma (x_b):\ t\ge 0\}.\tag 7.6$$
Conversely, for any $x_b\in \pom , y\in \ell^*_{x_b}$, we have 
$$h\le u_0 \ \ \ \text{in}\ \ \Om , \tag7.7$$ 
by virtue of (7.3), (taking
 $x_0 = x_b, y_0 =  y_{0,b}$). This proves that $u_1$ is indeed
 a $c$-convex extension of $u_0$ to $\Om^\delta$.
 
 To proceed further, we let $y_b$ be the unique point in $\pom^*$,
 where $\ell^*_{x_b}$ intersects $\pom^*$. Since $\omega^*$ is
 also uniformly $c^*$-convex, $\ell^*_{x_b}$  only intersects $\p\omega^*$
 at the initial point $y_{0,b}$. Actually, only the uniform $c$-convexity
 of $u_0$ is needed to justify this. Henceforth we restrict $\ell^*_{x_b}$
 to the segment joining  $y_{0,b}$ to $y_b$. From our argument above, the
 mapping from $x_b$ to $y_b$ is onto $\pom^*$. We claim it is also 
 one-to-one. For suppose the c-function $h =  h_{x_b,y_b}$ meets
 $\pom$ at another point $x^\prime$. By increasing $h$, we 
 infer the existence of a further boundary point $x^{\prime\prime}$, 
 which is a saddle
 point for $u_0$, contradicting (7.3). It follows
 then that the mapping from $x_b$ to $y_b$ is a $C^3$ diffeomorphism
 from $\pom$ to $\pom^*$. Next, if $B_r$ is a sufficiently small
 tangent ball of $\Om$ at $x_b$, it will also be uniformly $c$-convex
 so again by (7.3), we obtain
 $$h_{x_b, y_b}(x) > h_{x_b, y}(x)
    \ \ \ \forall \ x\in B_r, \ y\in \ell^*_{x_b}. \tag 7.8$$
  and thus we have 
$$u_1=\max_{x_b\in\pom} \{u_0, h_{x_b,y_b}\},
                       \ \ \ x\in\Om^\delta.\tag7.9$$

To complete the proof of Lemma 7.1, we need to show
that for each $x\in \Om^\delta - \Om$, there exists a unique
$x_b\in \pom$ where the maximum in (7.8) is attained. For
this purpose, we invoke the $c$-transform of $u_1$, 
$$v_0(y) =\sup \{c(x, y)-u_1(x):\ x\in\Om^\delta\},
                                 \ \ \ y\in\Om^*,\tag 7.10$$
which extends the $c$-transform of $u_0$ in $\omega^*$.
Moreover, by (7.7), we see that for $y\in \ell^*_{x_b}$, the sup is
attained at $x_b$. Hence
$$v_0(y)=c(x_b, y) - u_0(x_b)\ \ \ \forall\ y\in \ell^*_{x_b}. \tag 7.11$$
One easily verifies that $v_0$ is smooth in $\bom^*-
\p\omega^*$. Using (7.11) and arguing as before, we infer
that for any point $x\in{\Om^\delta - \bom}$ , there exists a
unique point $y_b\in\pom^*$ such that
$$\align
h^*_{x, y_b}(y) &= c(x, y)-c(x, y_b)+v_0(y_b)\\
& \le v_0(y)\ \ \ \forall\ y\in\bom^*.\tag 7.12\\
\endalign$$
Moreover $x$ lies on the $c$-segment, $\ell_{y_b}$
given by
$$c_y(\ell_{y_b}, y_b)=\{Dv_0(y_b)+t\gamma^*(y_b):\ t\in [0, \ol\delta]\},\tag 7.13$$
where $\gamma^*$ denotes the unit outer normal to $\pom^*$
and $\ol\delta$ is a small constant. Note that $Dv_0(y_b) = c_y(x_b,y_b)$
From (7.12), we
see that the maximum in (7.9) is attained at $x_b, y_b$ so
$$u_1(x) = c(x, y_b)-c(x_b, y_b)+u_0(x_b)\ \ \ x\in \ell_{y_b},\tag 7.14$$
with  $\chi_{u_1} (\ell_{y_b} - \{x_b\}) = y_b$, $\chi_{u_1}(x_b) = \ell^*_{x_b}$.
From the obliqueness of $\ell_{y_b}$ on $\pom$, we have that
the mapping from $x\in \Om^r$ to $x_b$ is one-to-one for
sufficiently small $r$. This completes the proof of Lemma 7.1.

From Lemma 7.1, we see that the function $u_1$ is smooth in 
$\Om^\delta - \pom$. Furthermore, with $\delta$  sufficiently small,
$u_1$ will be tangentially uniformly convex on $\pom^r$, that is
$$[D_{ij}u_1 - c_{ij}(\cdot, Tu_1)]  \tau_i\tau_j \ge \lambda_0,\tag 7.15$$
where $\tau$ is the unit tangent vector on  $\pom^r$ and 
$ \lambda_0$ a positive constant. To take care of the normal direction, we 
modify $u_1$ in $\Om^\delta  - \Om$, by setting
$$u = u_1 + bd^2,\tag 7.16$$
where $b$ is a positive constant and $d$ denotes distance from $\Om$.
Again for $\delta$ sufficiently small, we infer that $u$ satisfies 
$$[D_{ij}u - c_{ij}(\cdot, Tu)]  \xi_i\xi_j \ge \lambda_0,\tag 7.16$$
in $\Om^\delta - \pom$ for a further positive constant  $\lambda_0$,
and any unit vector $\xi$.

We complete the proof of Lemma 5.1 by mollification.
Let $\rho\in C_0^\infty(B_1(0))$ be a mollifier, namely $\rho$ is
a smooth, nonnegative, and radially symmetric function supported
in the unit ball $B_1(0)$ such that the integral
$\int_{B_1(0)}\rho=1$. We show that a mollification of $u$, given
by
$$\align
u_\eps(x)& =\rho * u\tag 7.17\\
  & =\int_{\R^n}\eps^{-n}\rho(\frac {x-y}\eps)u(y)\, dy\\
  &=\int_{\R^n} \rho(y)u(x-\eps y)\,dy\\
                        \endalign $$
is uniformly c-convex in $\Om^{\delta/2}$, provided $\eps<\frac
12\delta$ is sufficiently small and $x\in \Om^{\delta/2}$ (so that
the value of $u$ outside $\Om^\delta$ is irrelevant). Note that
the image of the c-normal mapping of $u_\eps$ in $\Om^{\delta/2}$
is a smooth perturbation of $\Om^*$, and so is also uniformly
c-convex provided $\eps>0$ is sufficiently small.

It is easy to verify that
$$D u_{\eps} (x)= \int_{\R^n} \rho(y)D u(x-\eps y)dy,\tag 7.18 $$
$$\align
D^2 u_{\eps}(x) &=
  \int_{\R^n} \rho(y)D^2 u(x-\eps y)dy
     +\int_{\pom} \frac {1}{\eps^n}\rho(\frac {x-y}\eps)
                   \  \gamma\cdot(D^+u-D^-u)(y) \tag 7.19\\
 &\ge\int_{\R^n} \rho(y) D^2 u(x-\eps y)dy,\\
  \endalign $$
where 
$$\align
 &D^+u(y)=\lim_{y'\not\in \bom, y'\to y} Du(y'),\\
 &D^-u(y)=\lim_{y'\in \Om, y'\to y} Du(y').\\
 \endalign $$
Since, $\omega^*\Subset \Om^*$,  we have
$$D^+_\gamma u-D^-_\gamma u\ge C_0>0\ \ \text{on}\ \ \pom \tag 7.20$$
for some positive constant $C_0$.  
We divide $\Om^{\delta/2}$ into three parts:
$\Om^{\delta/2}=U_1\cup U_2\cup U_3$, where
$$\align
 & U_1=\{x\in\Om^{\delta/2}:\ \dist(x, \pom)\ge\eps\},\\
 & U_2=\{x\in\Om^{\delta/2}:\ \dist(x, \pom)\in ((1-\sigma)\eps, \eps)\},\\
 & U_3=\{x\in\Om^{\delta/2}:\ \dist(x, \pom)\le \eps'\},\\
 \endalign $$
where $\sigma\in (\frac 12, 1)$ is a constant close to 1.
Since $u$ is smooth, uniformly c-convex away from $\pom$, $u^\eps$
is obviously smooth, uniformly c-convex in $U_1$ provided $\eps$
is sufficiently small. By taking $\sigma>0$ sufficiently close to 1, 
for any point$x_0\in U_2$, $Du^\eps(x_0)$ is a small
 perturbation of $Du(x_0)$. By (7.19), we also see that the matrix
$$\{D^2 u^\eps(x_0) - A(x_0, Du^\eps(x_0))\}>0, \tag 7.21$$
namely $u^\eps$ is smooth, uniformly c-convex in $U_2$.

Finally we verify (7.21) in $U_3$. For any point $x_0\in
U_3$, we choose a coordinate system such that $x_0=(0, \cdots, 0,
x_{0, n})$, the origin $0\in\pom$ and $\pom$ is tangent
to $\{x_n=0\}$. To verify (25), we first consider a tangential
direction $\tau$, namely $\tau$ is a unit vector tangential to
$\pom$ at $0$. Without loss of generality we assume that $\tau=(1,
0, \cdots, 0)$. Then we need to prove that
$$D_{11}u^{\eps}(x_0)  -A_{11}(x_0, Du^\eps(x_0))>0. \tag 7.22$$
By our choice of coordinates, $D{_1}u^\eps(x_0)$ is a small
perturbation of $D_1u(x_0)$. Hence it suffices to verify that
$$D_{11}u^{\eps}(x_0)-A_{11}(x_0, {D_1}u(x_0), D'u^\eps(x_0))>0, \tag 7.23$$
where $D'u^\eps=(D_2u^\eps, \cdots, D{_n}u^\eps)$. By A3w, $A_{11}$ is
convex in $D'u^\eps$, whence it follows readily that
$$A_{11}(x_0, D{_1}u(x_0), D'u^\eps(x_0))
 \le \int_{\R^n} \frac {1}{\eps^n}\rho(\frac{x_0-y}\eps)
   A_{11}(x_0, D{_1}u(x_0), D'u(y)) dy. $$
Inequality (7.22) now follows from (7.19) and (7.23). 
 Note that the argument also applies to any
direction $\eta$ provided $\eta\cdot \gamma$ is sufficiently small.
Next we observe from the second integral in  (7.19) that (7.21)
holds in the normal direction $\gamma = e_n$. Furthermore,
$$D_{nn}u^{\eps}(x_0)-A_{nn}(x_0, Du^\eps(x_0))\ge K\tag 7.24$$
for some $K$ as large as we want, provided $\eps$ is sufficiently
small.
Now suppose the least eigenvalue of the matrix (7.21) is achieved in
direction $\xi$. We can decompose $\xi=c_1\tau+c_2e_n$. If
$c_2\ge c_0$ for some constant $c_0>0$,  then the matrix (7.21) in
direction $\xi$ is positive by (7.24). Otherwise the proof of (7.22)
applies and we also see that the matrix (7.21) in direction $\xi$
is positive. By appropriate adjustment of $\Om$, we complete
the proof of Lemma 5.1 and consequently also
Theorems 1.2 and 1.3.

To conclude this section we show that Lemma 5.1 may be proved
independently of  the arguments in Section 5 by direct construction of a 
uniformly $c$-convex function, $u_0$. To do this we let $y_0$ be a
 point in $\Om^*$ and $u_0$ be
the $c^*$- transform of the function 
$$\psi(y)=-(r^2-|y-y_0|^2)^{1/2},\tag 7.25$$
given by
$$u_0(x)=\sup\{c(x, y)-\psi(y),\ \ \ y\in B_r(y_0)\},$$
for sufficiently small $r >0$.
Then $u_0$ is a locally uniformly c-convex function defined in
some ball $B_R(0)$, with $R\to\infty$ as $r\to 0$, and the image of its
c-normal mapping,
$$\omega^*:=Tu_0(\Om)\subset B_r(y_0), $$
where $Tu_0$ is a diffeomorphism between $\Om$ and $\omega^*$.
As $Tu_0$ is defined on the ball $B_R(0)\Supset\Om$, $\omega^*$
is a smooth domain. Locally $u_0$ is a smooth perturbation of the
c-function
$$h_0(\cdot) = c(\cdot, y_0)+a_0, \tag 7.26$$
for some constant  $a_0$.

\vskip20pt

\centerline {\bf  8. Examples}

\bigskip

We repeat and expand somewhat the examples in [MTW], taking account that our cost
functions are the negatives of those there.

{\bf Example 1.}
$$c(x,y) = - \sqrt{1+ |x-y|^2 }  \tag 8.1$$
Here the vector field $Y$ and matrix $A$ are given by
$$\align
 Y(x,p) &= x + {p \over \sqrt{1 - |p|^2}}  , \tag 8.2\\
A(x,p) &= A(p) = - \big ( 1 - |p|^2 \big ) ^{1/2}
 \big ( I - p \otimes p \big ). \\
\endalign $$
The cost function satisfies condition A3. We remark that condition
A1 is only satisfied for $|p| < 1$ but this does not  prohibit
application of our results as the boundedness of target  domain
$\Omega^*$ ensures that $|Du| < 1$ for solutions of (1.7),
(1.14). More generally the conditions $p,q\in\R^n$ in A1 may be replaced
by $p,q\in$ some convex sub-domain.

\bigskip

{\bf Example 2.}
$$ c(x,y) = - \sqrt{1 - |x - y|^2} \tag 8.3$$
Here $c$ is only defined for $|x-y| \le 1$. The vector field $Y$
and matrix $A$  are given by
$$\align
 Y(x,p) &= - x + {p \over \sqrt{1 + |p|^2}} ,\tag 8.4\\
 A(x,p) &= A(p) = \big ( 1 + |p|^2 \big ) ^{1/2}
 \big ( I + p \otimes p \big ).\\
\endalign $$
The cost function satisfies condition A3. In order to directly apply our
results we need to assume $\Omega$ and $\Omega^*$ are strictly
contained in a ball of radius 1.

\bigskip

{\bf Example 3.} Let $f, g \in C^2(\Omega) , C^2(\Omega^*)$
respectively and
$$ c(x,y) = x \cdot y + f(x) g(y) \; . \tag 8.5 $$
If $|\triangledown f.\triangledown g| < 1$, then $c$ satisfies
A1, A2. If $f, g $ are convex, then $c$ satisfies A3w, while if
$f,g $ are uniformly convex, then $c$ satisfies A3. As indicated
in [MTW], the function (8.5) is equivalent to the square of the
distance between points on the graphs of $f$ and $g$.

\bigskip

{\bf Example 4.  \quad Power costs.}
$$ c(x,y) = \pm {1 \over m} |x - y | ^m , \quad  m
 \ne 0 , \quad log |x-y| \; , \; m = 0). \tag 8.6 $$
For $m \ne 1$ and $x \ne y$, when $m <1$, the vector fields $Y$
and matrices $A$ are given by
$$\align
 Y(x,p) &= x \pm |p|^{{2-m\over m-1}} p,\tag 8.7\\
 A(x,p) &= A(p) = \pm \big \{ |p|^{{m-2 \over m-1}} I
  + (m-2) |p|^{-{m \over m-1}} p \otimes p \big \} .\\
\endalign $$
The only cases for which condition A3w  is satisfied are $m=2
(\pm)$ and $- 2\le m <1$ (+ only ). For the latter,
condition A3 holds for $- 2 < m <1$. To apply our
results directly in the latter cases, we need to assume $\Omega$
and $\Omega^*$ are disjoint.

In [MTW] we also considered the cost function
$$c(x, y)=-(1+|x-y|^2)^{p/2}\tag 8.8$$
for $1\le p\le 2$, extending Example 1 to $p>1$. We point out here
that these functions only satisfy A3 under the restriction
$|x-y|^2<\frac 1{p-1}$. This condition was omitted in [MTW].

\bigskip

{\bf Example 5. \quad  Reflector antenna problem}

Corresponding results and examples may be obtained on other
manifolds such as the spheres $S^n$. Indeed the considerations
 in [MTW] stemmed from the treatment of the
reflector antenna problem by Wang in [W1], which may be
represented as an optimal transportation problem on the sphere
$S^n$ with cost function
$$ c(x,y) = log (1 - x \cdot y ) ,\tag 8.9 $$
which is simply the spherical analogue of the case $m=0$ in
Example 4 above. The corresponding vector field $Y$ is now given
by
$$Y(x, p)=x-\frac {2}{1+|p|^2}(x+p), \tag 8.10$$
where now $p$ belongs to the tangent space of $S^n$ at $x$, while
the matrix $A$ is given by
$$A=\frac 12 (|p|^2-1)g_0-p\times p, \tag 8.11$$
where $g_0$ denotes the metric on $S^n$. See [W1, W2,GW, MTW] for
more details. When the domains $\Omega$ and $\Omega^*$ have
disjoint closures, and spherically uniformly convex boundaries, we
obtain the global regularity of potentials.

 We will defer further
examination and extensions to intersecting domains and other cost
functions in a future work. We also point out here that Example 4
provides regularity for quadratic cost functions on spheres when
the points x and y are sufficiently close.

\vskip20pt

 \baselineskip=12pt
\parskip=1pt

\newpage

\Refs\widestnumber\key{ABC}

\item {[C1]}  Caffarelli, L.:
              The regularity of mappings with a convex potential,
              J. Amer. Math. Soc., 5(1992), 99-104.

\item {[C2]}  Caffarelli, L.:
              Boundary regularity of maps with convex potentials II.
              Ann. of Math. 144 (1996), no. 3, 453--496.
              
\item {[C3]}  Caffarelli, L.:
              Allocation maps with general cost functions, 
              Lecture Notes in Pure and Appl. Math., 177 (1996), 29-35.             

\item {[D]}   Delano\"e, Ph.:
              Classical solvability in dimension two of the second
              boundary value problem associated with
              the Monge-Amp\`ere operator,
              Ann. Inst. Henri Poincar\'e, Analyse Non Lin\'eaire,
              8(1991), 443-457.

\item {[GM]}  Gangbo, W.,  McCann, R.J.:
              The geometry of optimal transportation,
              Acta Math.,  177(1996), 113-161.

\item {[GT]}  Gilbarg, D., Trudinger, N.S.:
              Elliptic partial differential equations of second order,
              Second Edition, Springer, Berlin, 1983.

\item {[GW]}  Guan, P. and Wang, X.J.:
              On a Monge-Amp\`ere equation arising in geometric optics,
              J. Diff. Geom. 48(1998), 205--223.

\item {[L]}  Loeper, G.,
             On the regularity of maps solutions of optimal transportation
             problems, preprint (2006).

\item {[LT]} Lieberman, G. M. and Trudinger, N.S.,
             Nonlinear oblique boundary value problems for nonlinear
             elliptic equations,
             Trans. Amer. Math. Soc. 295(1986), 509--546.

\item {[LTU]} Lions, P.L.,  Trudinger, N.S.,  Urbas, J.:
              Neumann problem for equations of Monge-Amp\`ere
              type,
              Comm. Pure Appl. Math., 39(1986), 539-563.

\item {[MTW]} Ma, X.N., Trudinger, N.S., and Wang, X-J.,
              Regularity of potential functions of the optimal
              transportation problem,
              Arch. Rat. Mech. Anal., 177(2005), 151-183.

\item {[RR]} Rachev, S.T., Ruschendorff, L.:
             Mass transportation problems,
             Springer, Berlin, 1998.

\item {[T1]}  Trudinger, N.S.,
             On the Dirichlet problem for Hessian equations,
             Acta Math., 175(1995), 151-164.

\item {[T2]}  Trudinger, N.S.,
             Lectures on nonlinear elliptic equations of second order,
             Lectures in Math. Sci. 9, Univ. Tokyo, 1995.

\item {[T3]}  Trudinger, N.S., 
              Recent  developments in elliptic partial differential equations
              of  Monge-Amp\`ere type, ICM,Madrid, 3(2006), 291-302.  

\item {[TW]} Trudinger, N.S., and Wang, X-J.,
              On strict convexity and continuous differentiability of potential 
              functions in optimal transportation, preprint (2006) .

\item {[U1]}  Urbas, J.:
              On the second boundary value problem for equations of
              Monge-Amp\`ere type,
              J. Reine Angew. Math., 487(1997), 115-124.

\item {[U2]} Urbas, J.,
              Mass transfer problems,
              Lecture Notes, Univ. of Bonn, 1998.

\item {[U3]} Urbas, J.,
             Oblique boundary value problems for equations of
             Monge-Amp\`ere type,
             Calc. Var. PDE, 7(1998),19--39.

\item {[V]}  Villani, C.,
             Topics in optimal transportation,
             Amer. Math. Soc., 2003.

\item {[W1]}  Wang, X.J.:
              On the design of a reflector antenna,
              Inverse Problems, 12(1996), 351-375.

\item {[W2]}  Wang, X.J.:
             On the design of a reflector antenna II,
             Calc. Var. PDE,  20(2004), 329-341.

\endRefs

\enddocument

\end